\documentclass[11pt]{amsart}

\usepackage{graphicx,pinlabel}

\newtheorem{theorem}{Theorem}[section]
\newtheorem{proposition}[theorem]{Proposition}

\newtheorem{corollary}[theorem]{Corollary}
\newtheorem*{TCC}{Tunnel Classification Conjecture}
\newtheorem*{KTP}{Known Tunnel Phenomena}

\newcommand{\Lk}{\operatorname{Lk}}

\newcommand\selambda{\,\,\genfrac{}{}{0pt}{2}{\searrow}{\lambda}\,\,}
\newcommand\nelambda{\,\,\genfrac{}{}{0pt}{2}{\lambda}{\nearrow}\,\,}
\newcommand\serho{\,\,\genfrac{}{}{0pt}{2}{\searrow}{\rho}\,\,}
\newcommand\nerho{\,\,\genfrac{}{}{0pt}{2}{\rho}{\nearrow}\,\,}
\newcommand\setau{\,\,\genfrac{}{}{0pt}{2}{\searrow}{\tau}\,\,}
\newcommand\netau{\,\,\genfrac{}{}{0pt}{2}{\tau}{\nearrow}\,\,}

\begin{document}

\title[Iterated splitting]
{Iterated splitting and the\\ classification of knot tunnels}

\author{Sangbum Cho}
\address{Department of Mathematics Education\\
Hanyang University\\
Seoul 133-791\\
Korea}
\email{scho@hanyang.ac.kr}

\author{Darryl McCullough}
\address{Department of Mathematics\\
University of Oklahoma\\
Norman, Oklahoma 73019\\
USA}
\email{dmccullough@math.ou.edu}
\urladdr{www.math.ou.edu/$_{\widetilde{\phantom{n}}}$dmccullough/}
\thanks{The second author was supported in part by NSF grant DMS-0802424}

\subjclass{Primary 57M25}

\date{\today}

\keywords{knot, tunnel, (1,1), torus knot, regular, splitting, 2-bridge}

\begin{abstract} For a genus-1 1-bridge knot in $S^3$, that is, a
$(1,1)$-knot, a middle tunnel is a tunnel that is not an upper or lower
tunnel for some $(1,1)$-position. Most torus knots have a middle tunnel,
and non-torus-knot examples were obtained by Goda, Hayashi, and
Ishihara. In a previous paper, we generalized their construction and
calculated the slope invariants for the resulting examples. We give an
iterated version of the construction that produces many more examples, and
calculate their slope invariants. If one starts with the trivial knot, the
iterated constructions produce all the $2$-bridge knots, giving a new
calculation of the slope invariants of their tunnels. In the final section
we compile a list of the known possibilities for the set of tunnels of a
given tunnel number 1 knot.
\end{abstract}

\maketitle

\section*{Introduction}
\label{sec:intro}

Genus-$2$ Heegaard splittings of the exteriors of knots in $S^3$ have been
a topic of considerable interest for several decades. They form a class
large enough to exhibit rich and interesting geometric behavior, but
restricted enough to be tractable. Traditionally such splittings are
discussed with the language of knot tunnels, which we will use from now
on.

The article \cite{CMtree} developed two sets of invariants that together give a
complete classification of all tunnels of all tunnel number 1 knots. One is
a finite sequence of rational ``slope'' invariants, the other a finite
sequence of ``binary'' invariants. The latter is trivial exactly when the
tunnel is a $(1,1)$-tunnel, that is, a tunnel that arises as the ``upper"
or ``lower" tunnel of a genus-$1$ $1$-bridge position of the knot. In the
language of \cite{CMtree}, the $(1,1)$-tunnels are called semisimple, apart from
those which occur as the well-known upper and lower tunnels of a $2$-bridge
knot, which are distinguished by the term ``simple". The tunnels which are
not $(1,1)$-tunnels are called regular.

For quite a long time, the only known examples of knots having both regular
and $(1,1)$-tunnels were (most) torus knots, whose tunnels were classified
by M. Boileau, M. Rost, and H. Zieschang~\cite{B-R-Z} and independently by
Y. Moriah~\cite{Moriah}. Recently, another example was found by H. Goda and
C. Hayashi~\cite{Goda-Hayashi}. The knot is the Morimoto-Sakuma-Yokota
$(5,7,2)$-knot, and Goda and Hayashi credit H.\ Song with bringing it to
their attention. Using his algorithm to compute tunnel invariants,
K. Ishihara verified that the tunnel is regular, and in view of this, we
refer to this example as the Goda-Hayashi-Ishihara tunnel. As noted
in~\cite{Goda-Hayashi}, a simple modification of their construction,
varying a nonzero integer parameter $n$, produces an infinite collection of
very similar examples.

In \cite{CMsplitting}, we gave an extensive generalization of the
Goda-Hayashi-Ishihara example, called the splitting construction, to
produce all examples directly obtainable by the geometric phenomenon that
underlies it. In addition, we gave an effective method to compute the full
set of invariants of the examples. Our construction will be reviewed in
Section~\ref{sec:splitting}.

In this paper, we develop an iterative method that begins with the result
of a splitting construction. The steps are not exactly splittings in the
sense of \cite{CMsplitting}, but are similar enough that we may call this
iterated splitting. The steps may be repeated an arbitrary number of times,
giving an immense collection of new examples of regular tunnels of
$(1,1)$-knots. At each step, a choice of nonzero integer parameter allows
further variation. Starting from each of the four splitting constructions,
we find two distinct ways to iterate, giving eight types of iteration.
Section~\ref{sec:iterated} describes the constructions in detail.

As with the splitting construction, the binary invariants of these new
tunnels are easy to find, but the slope invariants require more
effort. Fortunately, the general method given in \cite{CMsplitting} for
tunnels obtained by splitting can be applied to obtain the slope invariants
for the iterated construction, as we detail in
Section~\ref{sec:iterated_slopes}. The method is effective and could easily
be programmed to read off slope invariants at will.

The iterated splitting construction actually sits in plain view in a very
familiar family of examples, the semisimple tunnels of $2$-bridge knots. In
Section~\ref{sec:2bridge_iteration}, we present a special case of the
iterated splitting method that, as one varies its parameters, produces all
semisimple tunnels of all $2$-bridge knots. No doubt there is a geometric
way to verify this, but our proof is short and entirely algebraic: we
simply calculate the slope sequences of the tunnels produced by the
iterations and see that they are exactly the sequences that arise from this
class of tunnels. The binary invariants are trivial in both cases, and
since the invariants together form a complete invariant of a knot tunnel,
the verification is complete.

The work in this paper greatly enlarges the list of known examples of
tunnels having a pair of $(1,1)$-tunnels and an additional regular tunnel,
motivating us to compile a list of known phenomena for the set of tunnels
of a given tunnel number~$1$ knot. In the final section, we give the list
of seven known cases, which includes three new cases apparent from examples
recently found by John Berge using his software package \textit{Heegaard.}
The authors are very grateful to John, not only for the new examples, but
also for providing patient consultation to help us understand his methods.

Although we do provide a review of the splitting construction of
\cite{CMsplitting}, this paper presupposes a reasonable familiarity with
that work. We have not included a review of the general theory
of~\cite{CMtree}, as condensed reviews are already available in several of
our articles. For the present paper, we surmise that Section~1
of~\cite{CMgiant_steps} together with the review sections
of~\cite{CMsemisimple} form the best option for most readers.

\section{The splitting construction}
\label{sec:splitting}

\begin{figure}
\begin{center}
\labellist
\pinlabel $m$ at 152 55
\pinlabel $\ell$ at 190 55
\pinlabel \Large $T$ at 0 140
\small\pinlabel $T_{3,5}$ at 212 27
\endlabellist
\includegraphics[width=0.5\textwidth]{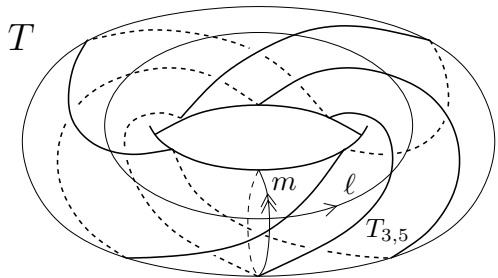}
\caption{$\ell$, $m$, and $T_{3,5}$.}
\label{fig:torus_knot_notation}
\end{center}
\end{figure}

In this section we will review the splitting construction
from~\cite{CMsplitting}. To set notation,
Figure~\ref{fig:torus_knot_notation} shows a standard Heegaard torus $T$ in
$S^3$, and an oriented longitude-meridian pair $\{\ell,m\}$ which will be
our ordered basis for $H_1(T)$ and for the homology of a product
neighborhood $T\times I$. For a relatively prime pair of integers $(p,q)$,
we denote by $T_{p,q}$ a torus knot isotopic to a $(p,q)$-curve in $T$. In
particular, $\ell=T_{1,0}$ and $m=T_{0,1}$, also $T_{p,q}$ is isotopic in
$S^3$ to $T_{q,p}$ in $S^3$, $T_{-p,-q}=T_{p,q}$ since our knots are
unoriented.

\begin{figure}
\begin{center}
\labellist
\pinlabel \large (a) at 113 3
\pinlabel $\rho$ at 163 286
\pinlabel $\lambda$ at 28 185
\pinlabel $\tau$ at 91 143
\pinlabel $\sigma$ at 48 194
\pinlabel $\lambda$ at 172 114
\pinlabel $\sigma$ at 152 98
\pinlabel $\rho$ at 37 12
\pinlabel $K_\tau$ at -11 133
\pinlabel $K_\lambda$ at -11 37
\pinlabel $K_\lambda$ at 192 254
\pinlabel $K_\tau$ at 192 156
\pinlabel \large (b) at 337 3
\pinlabel $\lambda$ at 267 212
\pinlabel $\rho$ at 334 212
\pinlabel $\lambda$ at 400 212
\pinlabel $K_\tau$ at 441 189
\pinlabel $K_\rho$ at 441 167
\pinlabel $K_\rho$ at 227 100
\pinlabel $K_\lambda$ at 441 78
\pinlabel $\sigma$ at 279 133
\pinlabel $\tau$ at 267 56
\pinlabel $\tau$ at 401 56
\endlabellist
\includegraphics[width=0.8\textwidth]{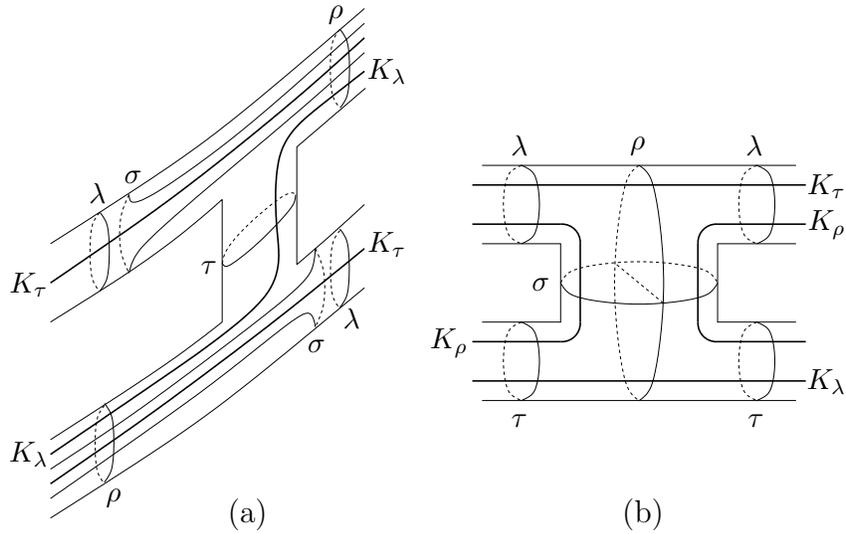}
\caption{The drop-$\lambda$ disk $\sigma$, first as seen in a
neighborhood of $K_\tau=T_{p+r,q+s}$ and the tunnel $\tau$, then after
dropping $K_\lambda=T_{r,s}$ and part of $K_\rho=T_{p,q}$.}
\label{fig:drop-lambda}
\end{center}
\end{figure}

Four kinds of disks, called drop-$\lambda$, lift-$\lambda$, drop-$\rho$,
and lift-$\rho$ disks, are used in the splitting construction.
Figure~\ref{fig:drop-lambda}(a) shows a torus knot $T_{p+r,q+s}$, its
middle tunnel $\tau$, the principal pair $\{\lambda,\rho\}$ of $\tau$, the
knots $K_\rho=T_{p,q}$, and $K_\lambda=T_{r,s}$, and a drop-$\lambda$
disk, called $\sigma$ there. Figure~\ref{fig:drop-lambda}(b) is an isotopic
repositioning of the configuration of Figure~\ref{fig:drop-lambda}(a): the
vertical coordinate is the $I$-coordinate in a product neighborhood
$T\times I$, $K_\tau$ and $K_\lambda$ lie on concentric tori in $T\times
I$, and the $1$-handle with cocore $\sigma$ is a vertical $1$-handle
connecting tubular neighborhoods of these two knots. The term
``drop-$\lambda$'' is short for ``drop-$K_\lambda$'', motivated by the fact
that a copy of $K_\lambda$ can be dropped to a lower torus level, as in
Figure~\ref{fig:drop-lambda}(b).

\begin{figure}
\begin{center}
\labellist

\pinlabel \large (a) at 113 3
\pinlabel $\rho$ at 163 286
\pinlabel $\lambda$ at 28 185
\pinlabel $\tau$ at 91 143
\pinlabel $\sigma$ at 48 194
\pinlabel $\lambda$ at 172 114
\pinlabel $\sigma$ at 152 98
\pinlabel $\rho$ at 37 12
\pinlabel $K_\tau$ at -11 133
\pinlabel $K_\lambda$ at -11 37
\pinlabel $K_\lambda$ at 192 254
\pinlabel $K_\tau$ at 192 156
\pinlabel \large (b) at 337 3
\pinlabel $\tau$ at 267 212
\pinlabel $\rho$ at 334 212
\pinlabel $\tau$ at 400 212
\pinlabel $K_\lambda$ at 441 189
\pinlabel $K_\rho$ at 441 167
\pinlabel $K_\rho$ at 227 100
\pinlabel $K_\tau$ at 441 78
\pinlabel $\sigma$ at 283 133
\pinlabel $\lambda$ at 267 56
\pinlabel $\lambda$ at 401 56
\endlabellist
\includegraphics[width=0.8\textwidth]{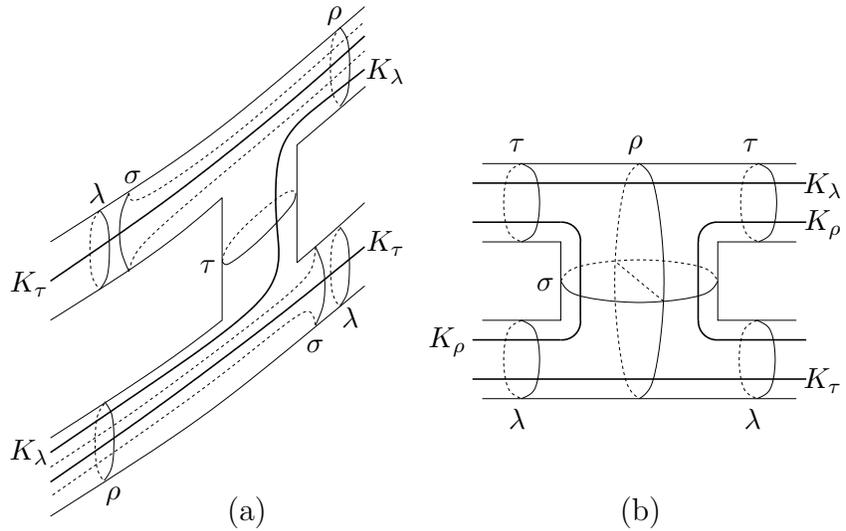}
\caption{The lift-$\lambda$ disk $\sigma$, first as seen in a
neighborhood of $K_\tau=T_{p+r,q+s}$ and the tunnel $\tau$, then after
lifting $K_\lambda=T_{r,s}$ and part of $K_\rho=T_{p,q}$.}
\label{fig:lift-lambda}
\end{center}
\end{figure}

A lift-$\lambda$ disk is similar, and is shown in
Figure~\ref{fig:lift-lambda}. Drop-$\rho$ and lift-$\rho$ disks are
similar, except that they cut across the upper copy of $\lambda$, travel
over the portion of the neighborhood of $T_{p+r,q+s}$ that does not contain
the drop-$\lambda$ disks, and cut across the lower copy of $\lambda$, while
staying disjoint from the copies of~$\rho$.

\begin{figure}
\begin{center}
\labellist
\pinlabel $\lambda$ at 24 119
\pinlabel $\gamma_n$ at 80 119
\pinlabel $\lambda$ at 133 119
\pinlabel $K_\tau=T_{p+r,q+s}$ at 200 101
\pinlabel $K_\rho=T_{p,q}$ at 187 80
\pinlabel $\tau$ at 24 -7
\pinlabel $\tau$ at 134 -7
\pinlabel $K_\lambda=T_{r,s}$ at 187 10
\endlabellist
\includegraphics[width=0.36\textwidth]{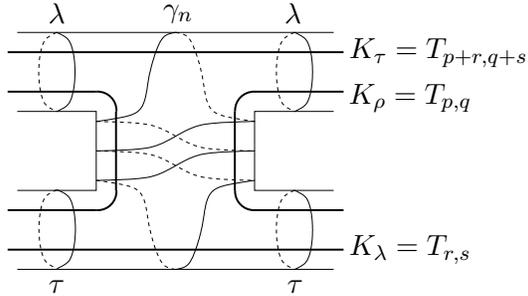}
\caption{The disk $\gamma_n$ is obtained from $\rho$ by $n$ right-handed
half-twists along $\sigma$. The case $n=3$ is shown here. For $n<0$, the
half-twists are left-handed, while $\gamma_0=\rho$.}
\label{fig:gamma}
\end{center}
\end{figure}

The splitting constructions split off a copy of $K_\rho=T_{p,q}$ or
$K_\lambda=T_{r,s}$ from $K_\tau=T_{p+r,q+s}$, producing copies of these
knots on two concentric torus levels, then sum the copies together by a
pair of arcs with some number of twists. In the case of the the
drop-$\lambda$ splitting, the first step was illustrated in
Figure~\ref{fig:drop-lambda}. Next, consider the disk $\gamma_n$ shown in
Figure~\ref{fig:gamma}. It is obtained from $\rho$ by $n$ right-handed
half-twists along $\sigma$. When $n<0$, the twists are left-handed, while
$\gamma_0=\rho$. The $\gamma_n$ are nonseparating, since each meets
$K_\tau$ in a single point.

Each $\gamma_n$ with $n\neq 0$ is a tunnel for the knot obtained by joining
the copies of $K_\tau$ and $K_\lambda$ in Figure~\ref{fig:drop-lambda} by a
pair of vertical arcs that have $n$ right-handed half-twists. That is, for
$n\neq 0$ going from $\tau$ to $\gamma_n$ is a cabling construction
replacing $\rho$, so that the principal pair of $\gamma_n$ is
$\{\lambda,\tau\}$. The case of $n=0$ does not produce a cabling
construction (that is, the resulting tunnel would be $\rho$ so the
principal path would have reversed direction).

The lift-$\lambda$, drop-$\rho$, and lift-$\rho$ splittings are exactly
analogous, using the lift-$\lambda$, drop-$\rho$, and lift-$\rho$ disks
as $\sigma$ in the respective cases.

The slope invariants of the resulting tunnels are the slopes of the disks
$\gamma_n$ in certain coordinates. To calculate them, we need the slopes of
the drop- and lift-disks. We review the method used in~\cite{CMsplitting},
which will apply to the iterated construction that we will develop in this
paper.

\begin{figure}
\begin{center}
\labellist
\pinlabel \large (a) at 67 -7
\pinlabel \large (b) at 217 -7
\pinlabel \large (c) at 370 -7
\scriptsize\pinlabel $D_U^-$ at 20 123
\scriptsize\pinlabel $D$ at 65 123
\scriptsize\pinlabel $D_U^+$ at 110 123
\scriptsize\pinlabel $D_U^-$ at 170 123
\scriptsize\pinlabel $D_U^+$ at 261 123
\scriptsize\pinlabel $D_U^-$ at 320 123
\scriptsize\pinlabel $D^0$ at 368 123
\scriptsize\pinlabel $D_U^+$ at 418 123
\scriptsize\pinlabel $K_U$ at -8 105
\scriptsize\pinlabel $K_U^0$ at 142 95
\scriptsize\pinlabel $K_L^0$ at 142 39
\scriptsize\pinlabel $K_L$ at -8 31
\small\pinlabel $\sigma$ at 28 68
\scriptsize\pinlabel $D_L^-$ at 20 12
\scriptsize\pinlabel $D_L^+$ at 111 12
\scriptsize\pinlabel $D_L^-$ at 171 12
\scriptsize\pinlabel $D_L^+$ at 262 12
\scriptsize\pinlabel $D_L^-$ at 322 12
\scriptsize\pinlabel $D_L^+$ at 419 12
\endlabellist
\includegraphics[width=0.93\textwidth]{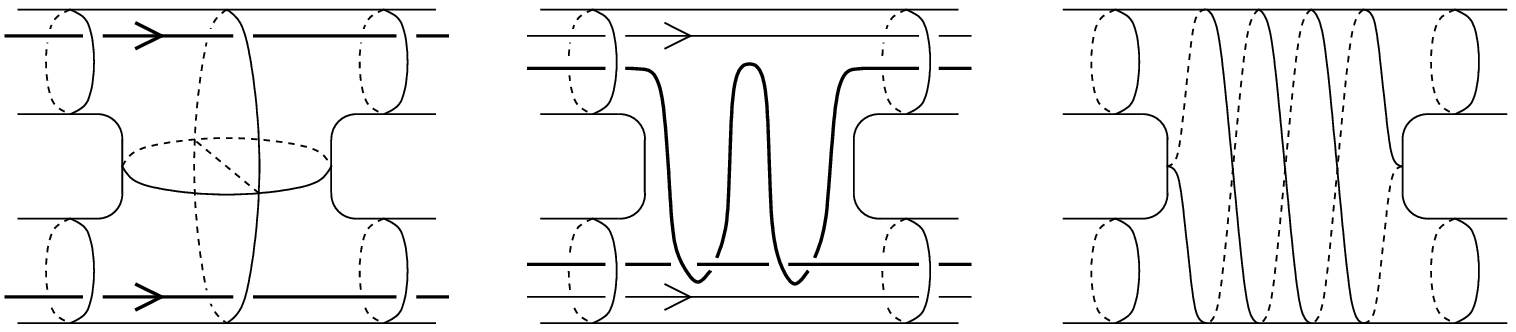}
\caption{The setup for the first general slope calculation.}
\label{fig:first_general}
\end{center}
\end{figure}

Figure~\ref{fig:first_general} illustrates the setup for the slope
calculation. The first drawing shows tubular neighborhoods of two
(oriented) knots $K_U$ and $K_L$, contained in a product neighborhood
$T\times I$ of a Heegaard torus $T$ of $S^3$. The neighborhoods are
connected by a vertical $1$-handle to yield a genus-$2$ handlebody $H$. In
our context, $H$ will always be unknotted, although that is not needed for
the calculations of this and the next section.

We interpret $K_U$ as the ``upper'' knot, contained in $T\times [0,1/4)$
and $K_L$ as the ``lower'' knot, contained in $T\times (3/4,1]$ (the
$I$-coordinate of $T\times I$ increases as one moves downward in our
figures). The vertical $1$-handle with cocore $\sigma$ is assumed to run
between $T\times \{1/4\}$ and $T\times \{3/4\}$, with $\sigma$ as its
intersection with~$T\times \{1/2\}$.

The homology group $H_1(T\times I)\cong H_1(T)$ will have ordered basis the
oriented longitude and meridian $\ell$ and $m$ shown in
Figure~\ref{fig:torus_knot_notation}. Our linking convention is that
$\Lk(m\times \{1\},\ell\times \{0\})=+1$. Now, suppose that $K_U$ represents
$(\ell_U,m_U)$ and $K_L$ represents $(\ell_L,m_L)$ in $H_1(T\times I)$.
Since $\Lk(m\times \{0\},\ell\times \{1\})=0$, we have
$\Lk(K_U,K_L)=m_U\ell_L$.

The disks $D_U^+$ and $D_U^-$ in Figure~\ref{fig:first_general} are
parallel in $H$, as are the disks $D_L^+$ and $D_L^-$, and these four disks
bound a ball $B$. Figure~\ref{fig:first_general}(a) shows a slope disk
$D$. Associated to $D$ is a slope-$0$ separating disk $D^0$, defined by the
requirement that it meets $D$ in a single arc and the core circles of its
complementary solid tori in $H$ have linking number $0$ in $S^3$. For this
setup, \cite[Proposition 5.1]{CMsplitting} tells us the slope $m_\sigma$ of $\sigma$
in $(D,D^0)$-coordinates.
\begin{proposition}\label{prop:first_general_slope_calculation}
In Figure~\ref{fig:first_general}, the slope $m_\sigma$ of $\sigma$ in
$(D,D^0)$-coordinates is $2\Lk(K_U,K_L)$. Consequently, if $K_U$ represents
$(\ell_U,m_U)$ and $K_L$ represents $(\ell_L,m_L)$ in $H_1(T\times I)$,
then $m_\sigma$ equals $2m_U\ell_L$.
\end{proposition}

Proposition~6.1 of \cite{CMsplitting} then gives the slope of $\gamma_n$.
\begin{proposition}\label{prop:second_general_slope_calculation}
The slope of $\gamma_n$ in $(D,D^0)$-coordinates is $m_\sigma+1/n$.
\end{proposition}

As detailed in \cite[Proposition 7.1]{CMsplitting}, applying
Proposition~\ref{prop:first_general_slope_calculation}
to splitting disks gives their slopes
in terms of $p$, $q$, $r$, and $s$.
\begin{corollary}\label{coro:splitting_disk_slopes}
The slopes of the splitting disks are as follows:
\begin{enumerate}
\item[(a)] In $(\rho,\rho^0)$-coordinates,
the drop-$\lambda$ disk has slope $2r(q+s)$.
\item[(b)] In $(\rho,\rho^0)$-coordinates,
the lift-$\lambda$ disk has slope $2s(p+r)$.
\item[(c)] In $(\lambda,\lambda^0)$-coordinates,
the drop-$\rho$ disk has slope $2p(q+s)$
\item[(d)] In $(\lambda,\lambda^0)$-coordinates,
the lift-$\rho$ disk has slope $2q(p+r)$.
\end{enumerate}
\end{corollary}
Proposition~\ref{prop:second_general_slope_calculation} then gives immediately
the slopes of the tunnels obtained by splitting constructions
using $\gamma_n$.
\begin{proposition}\label{prop:sigma_slopes}
For the torus knot $T_{p+r,q+s}$:
\begin{enumerate}
\item[(a)] A drop-$\lambda$ splitting has slope $2r(q+s)+1/n$.
\item[(b)] A lift-$\lambda$ splitting has slope $2s(p+r)+1/n$.
\item[(c)] A drop-$\rho$ splitting has slope $2p(q+s)+1/n$.
\item[(d)] A lift-$\rho$ splitting has slope $2q(p+r)+1/n$.
\end{enumerate}
\end{proposition}

\section{The iterated splitting construction}
\label{sec:iterated}

\begin{figure}
\begin{center}
\labellist
\pinlabel \large (a) at 107 413
\small\pinlabel $\rho$ at 166 583
\small\pinlabel $K_\tau$ at 205 565
\small\pinlabel $K_{\gamma_{n_0}^0}$ at 212 480
\small\pinlabel $\gamma_{n_0}^0$ at 105 442
\small\pinlabel $\tau$ at 167 443
\small\pinlabel $K_\rho$ at 205 460
\normalsize\pinlabel \large (b) at 107 217
\small\pinlabel $\rho$ at 166 394
\small\pinlabel $K_{\gamma_{n_0}^0}$ at 212 291
\small\pinlabel $K_\rho$ at 207 269
\small\pinlabel $\gamma_{n_0}^0$ at 105 251
\small\pinlabel $\tau$ at 167 253
\normalsize\pinlabel \large (c) at 107 8
\small\pinlabel $\rho$ at 166 192
\small\pinlabel $K_{\gamma_{n_0}^0}$ at 232 107
\small\pinlabel $K_\rho$ at 227 52
\small\pinlabel $\gamma_{n_0}^0$ at 87 30
\small\pinlabel $\tau$ at 164 32
\normalsize\pinlabel \large (d) at 398 359
\small\pinlabel $\rho$ at 452 554
\small\pinlabel $K_{\gamma_{n_1}^1}$ at 520 410
\small\pinlabel $\gamma_{n_0}^0$ at 374 392
\small\pinlabel $\gamma_{n_1}^1$ at 451 392
\normalsize\pinlabel \large (e) at 398 60
\small\pinlabel $\rho$ at 439 293
\small\pinlabel $K_{\gamma_{n_1}^1}$ at 537 147
\small\pinlabel $K_\rho$ at 534 120
\small\pinlabel $\gamma_{n_0}^0$ at 346 101
\small\pinlabel $\gamma_{n_1}^1$ at 440 101
\endlabellist
\includegraphics[width=0.9\textwidth]{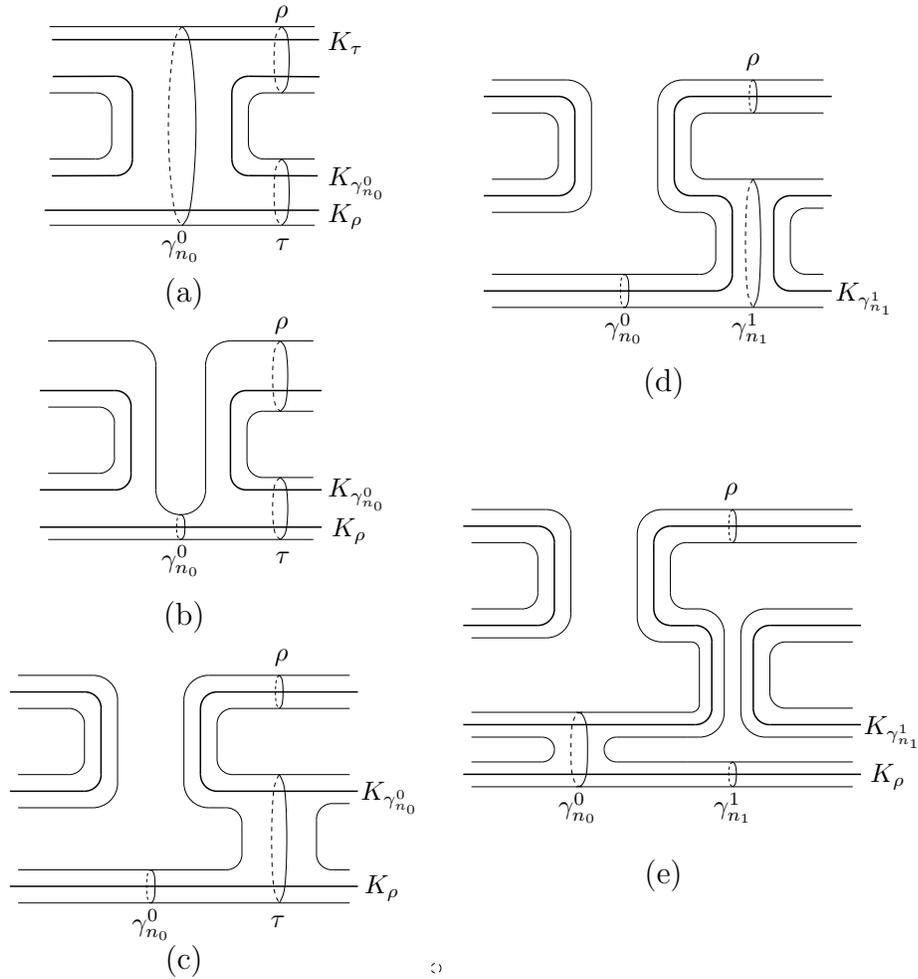}
\caption{The first case of the drop-$\rho$ iteration.}
\label{fig:iterate}
\end{center}
\end{figure}

We are now prepared to describe the iterated splitting construction.
We begin with the drop-$\rho$ case, as it is the case we will need in our
later application to $2$-bridge knots in Section~\ref{sec:2bridge_iteration}.
Figure~\ref{fig:iterate}(a) shows a
knot resulting from a drop-$\rho$ splitting. Its tunnel will now be denoted
by $\gamma_{n_0}^0$, the superscript distinguishing it from later
tunnels. Its principal pair $\{\rho,\tau\}$ is also shown.

In $S^3$, $\gamma_{n_0}^0$ would appear with twists along the horizontal
drop-$\rho$ disk $\sigma$, so Figure~\ref{fig:iterate}(a) is only a
picture up to abstract homeomorphism. Nonetheless, the vertical coordinate
represents the levels of $T\times I$, which
will be true in the remaining drawings of Figure~\ref{fig:iterate}, so it
will be seen that knots in $1$-bridge position will always be obtained.

In Figure~\ref{fig:iterate}(b), $\gamma_{n_0}^0$ and a portion of the
surrounding handlebody $H$ have been shrunk vertically, keeping
$K_{\gamma_{n_0}^0}$ fixed. The horizontal line at the bottom is a copy of
$K_\rho$, as indicated. The picture of $\gamma_{n_0}^0$ without twisting is
now accurate, but in the true picture in $S^3$, the two vertical
$1$-handles would be intertwined by $n_0$ right-hand half-twists rather
than being straight. The bottom part of the picture in $S^3$, from the
level $K_{\gamma_{n_0}^0}$ and below, is as seen in
Figure~\ref{fig:iterate}(b).

Figure~\ref{fig:iterate}(c) is obtained from Figure~\ref{fig:iterate}(b) by
an isotopy of $H$, keeping $K_{\gamma_{n_0}^0}$ and $K_\rho$ fixed. The
effect is to create the setup picture of Figure~\ref{fig:first_general}(a)
near $\tau$, with $K_U= K_{\gamma_{n_0}^0}$ and $K_L=K_\rho$. Notice that
in the orientations needed for the first general slope calculation,
$K_\rho$ is oriented left-to-right, and $K_{\gamma_{n_0}^0}$ must be
oriented so that the portion that intersects $\tau$ and originally came
from the copy of $K_\rho$ in the splitting construction used to create
$K_{\gamma_{n_0}^0}$ is also oriented from left-to-right. With this
orientation on $K_{\gamma_{n_0}^0}$ the top portion that originally came
from $K_\tau$ will be oriented from left-to-right or from right-to-left
according as $n_0$ is odd or even. This will be a key observation when we
compute the slope invariants of the iterated splitting constructions in
Sections~\ref{sec:iterated_slopes}.

Figure~\ref{fig:iterate}(d) differs from Figure~\ref{fig:iterate}(c) only
in that $\tau$ has been replaced by $\gamma_{n_1}^1$, which in $S^3$ would be
seen with $n_1$ right-hand half-twists. This is a cabling construction.
The resulting knot $K_{\gamma_{n_1}^1}$ is in $1$-bridge position, and was
obtained from $K_{\gamma_{n_0}^0}$ and the copy of $K_\rho$ by connecting
them with two vertical arcs with $n_1$ half-twists. The principal pair of
$\gamma_{n_1}^1$ is~$\{\rho,\gamma_{n_0}^0\}$.

The stage is now set to repeat the construction using $\gamma_{n_0}^0$ and
$\gamma_{n_1}^1$ in the role of $\tau$ and $\gamma_{n_0}^0$ in the previous
step. Figure~\ref{fig:iterate}(e) is obtained from
Figure~\ref{fig:iterate}(d) two steps, analogous to the steps from
Figure~\ref{fig:iterate}(a) to Figure~\ref{fig:iterate}(b) and from
Figure~\ref{fig:iterate}(b) to Figure~\ref{fig:iterate}(c).  First,
$\gamma_{n_1}^1$ is shrunk vertically, then $H$ is moved as indicated,
creating the setup picture of Figure~\ref{fig:first_general}(a) in the
lower left-hand area of Figure~\ref{fig:iterate}(d). Again, in $S^3$ the
two vertical $1$-handles in the middle would be intertwined with $n_1$
half-twists. Another copy of $K_\rho$ appears at the bottom.

The next cabling construction replaces $\gamma_{n_0}^0$ by
$\gamma_{n_2}^2$, and $K_{\gamma_{n_2}^2}$ is obtained by joining
$K_{\gamma_{n_1}^1}$ and the copy of $K_\rho$ with two vertical arcs with
$n_2$ half-twists. The principal pair of $\gamma_{n_2}^2$ is
$\{\rho,\gamma_{n_1}^1\}$. The true picture in $S^3$ has $n_0$ half-twists
in the two vertical $1$-handles connecting the top and second levels of
Figure~\ref{fig:first_general}(e), $n_1$ half-twists in the two
vertical $1$-handles connecting the second and third levels, and
$\gamma_{n_2}^2$ appears with $n_2$ half-twists.

The iteration can be continued indefinitely, producing a sequence of
tunnels $\gamma_{n_m}^m$ with principal pairs
$\{\rho,\gamma_{n_{m-1}}^{m-1}\}$, and the knots $K_{\gamma_{n_m}^m}$ in
$(1,1)$-position.

We indicate this sequence by $\tau \serho \gamma_{n_0}^0 \serho
\gamma_{n_1}^2 \serho \cdots$. The cabling constructions in the iterations
all retain $\rho$ in their principal pairs so have binary invariant $0$,
although the original drop-$\rho$ splitting that produces $\gamma_{n_0}^0$
may have nontrivial binary invariant.

\begin{figure}
\begin{center}
\labellist
\small\pinlabel $\rho$ at 185 164
\small\pinlabel $\gamma_{n_0}^0$ at 109 169
\small\pinlabel $K_\tau$ at 236 142
\small\pinlabel $K_{\gamma_{n_0}^0}$ at 245 110
\pinlabel $\tau$ at 185 -9
\endlabellist
\includegraphics[width=0.4\textwidth]{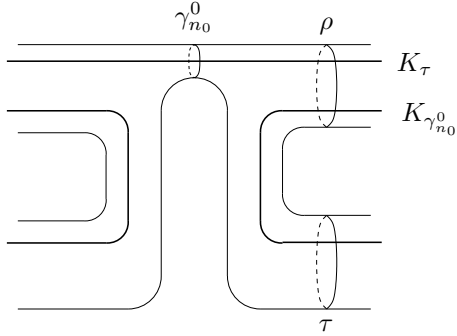}
\caption{The second case of the drop-$\rho$ iteration.}
\label{fig:tau_iterate}
\end{center}
\end{figure}

From Figure~\ref{fig:iterate}(a) there is a second way to
proceed. Figure~\ref{fig:tau_iterate} shows an alternative to the isotopy
in Figure~\ref{fig:iterate}(b), that shrinks $\gamma_{n_0}^0$ upward. The
next step replaces $\rho$ by $\gamma_{n_1}^1$, which has principal pair
$\{\tau,\gamma_{n_0}^0\}$, and $K_{\gamma_{n_1}^1}$ is obtained by joining
copies of $K_{\gamma_{n_0}^0}$ and $K_\tau$ by vertical arcs. The
successive iterations each add on another copy of $K_\tau$, moving upward,
and retain $\tau$ in their principal pairs. We indicate this sequence by
$\tau \serho \gamma_{n_0}^0 \netau \gamma_{n_1}^1 \netau \gamma_{n_2}^2
\netau \cdots$. The up-or-down direction of the diagonal arrow indicates
whether the knot that is joined to the previous one is a copy of the
original $K_U$ (in this case, $K_\tau$) or the original $K_L$ (in this
case, $K_\rho$), and the letter above it indicates which of $\rho$,
$\lambda$, or $\tau$ is retained in the principal pair.

Starting with the drop-$\lambda$ splitting instead of the drop-$\rho$
splitting produces two more interations,
\[\tau \selambda \gamma_{n_1}^1 \selambda
\gamma_{n_2}^2 \selambda\cdots\text{\ and\ }
\tau \selambda \gamma_{n_1}^1 \netau
\gamma_{n_2}^2 \netau \gamma_{n_3}^3 \netau \cdots\ .\]

Starting with the lift-$\rho$ splitting instead of the drop-$\rho$
splitting produces two more,
\[\tau \nerho \gamma_{n_1}^1 \nerho
\gamma_{n_2}^2 \nerho\cdots\text{\ and\ } \tau \nerho \gamma_{n_1}^1 \setau
\gamma_{n_2}^2 \setau \gamma_{n_3}^3 \setau \cdots\ ,\] and starting these
with the lift-$\lambda$ splitting give the latter two but with $\lambda$
replacing $\rho$. Provided that one started with a tunnel $\tau$
which was not trivial and not simple, the eight sequences are distinct,
since they have distinct principal paths.

\section{The iterated splitting slope invariants}
\label{sec:iterated_slopes}

We begin with the slope invariants. Consider the first iteration discussed
in Section~\ref{sec:iterated}, whose initial steps were illustrated in
Figure~\ref{fig:iterate}. The initial step is a regular drop-$\rho$
splitting, and according to Proposition~\ref{prop:sigma_slopes}(c), the
slope of the resulting tunnel disk $\gamma_{n_0}^0$ is $2p(q+s)+1/n_0$.

The first iterate $\gamma_{n_1}^1$ is obtained using the setup of
Figure~\ref{fig:first_general}(a) with $K_U=K_{\gamma_{n_0}^0}$ and
$K_L=K_\rho$ as in Figure~\ref{fig:iterate}(c). Now $K_{\gamma_{n_0}^0}$ is
obtained by connecting $K_{\tau}=T_{p+r,q+s}$ and $K_\rho=T_{p,q}$ with two
arcs intertwined with $n_0$ half-twists. The portion of
$K_{\gamma_{n_0}^0}$ seen in setup picture for calculating the slope of
$\gamma_{n_1}^1$ must be oriented from left-to-right, so it is obtained by
adding the left-to-right orientation of $K_\rho$ to either the
left-to-right or right-to-left orientation of $K_\tau$, according as $n_0$
is odd or even. In $H_1(T\times I)$, $K_\tau$ (with left-to-right
orientation) represents $(p+r,q+s)$ and $K_\rho$ represents $(p,q)$, so
$K_{\gamma_{n_0}^0}$ with this orientation represents
$(p,q)+(-1)^{1+n_0}(p+r,q+s) =(1+(-1)^{1+n_0})(p,q) +
(-1)^{1+n_0}(r,s)$. Therefore $\Lk(K_{\gamma_{n_0}^0},K_\rho)
=p((1+(-1)^{1+n_0}q + (-1)^{1+n_0}s)$, and by
Proposition~\ref{prop:first_general_slope_calculation} the slope of
$\gamma_{n_1}^1$ in $(\tau,\tau^0)$-coordinates is $2pq(1+(-1)^{1+n_0}) +
2p s\, (-1)^{1+n_0} +1/n_1$.

To continue this process, let us put $\epsilon(k)=(-1)^{1+n_k}$,
$t(r,k)=\epsilon(r)\epsilon(r+1)\cdots\epsilon(k-1)$ for $r<k$, and
$t(k,k)=0$. Now, define
\[a(k)=t(0,k)\text{ and }A(k)=1+\sum_{r=0}^{k-1}t(r,k)\ .\]
In particular, $a(0)=A(0)=1$, $a(1)=\epsilon(0)$, $A(1)=1+\epsilon(0)$,
and since $\epsilon(k)t(r,k)=t(r,k+1)$,
\[\epsilon(k)a(k)=a(k+1)\text{ and }1+\epsilon(k)A(k)=A(k+1)\ .\]

We orient each $K_{\gamma_{n_k}^k}$ so that
the portion that came from $K_L=K_\rho$
is left-to-right, as this is the orientation needed
in order to compute the slope of $\gamma_{n_{k+1}}^{k+1}$ in the setup
of Figure~\ref{fig:first_general}(a).
In $H_1(T\times I)$, $K_{\gamma_{n_0}^0}$ represents
$(p,q)+\epsilon(0)(p+r,q+s)=A(1)(p,q)+a(1)(r,s)$.
For $k\geq 1$ assume inductively that
$K_{\gamma_{n_{k-1}}^{k-1}}$ represents $A(k)\,(p,q)+a(k)\,(r,s)$.
In the orientation on $K_{\gamma_{n_k}^k}$, the direction
on the portion from
$K_{\gamma_{n_{k-1}}^{k-1}}$ must be reversed exactly when $n_k$ is
even. Therefore in $H_1(T\times I)$,
$K_{\gamma_{n_k}^k}$ represents
\[(p,q)+\epsilon(k)(A(k)\,(p,q) + a(k)\,(r,s))
=A(k+1)\,(p,q)+a(k+1)\,(r,s)\ ,\]
completing the induction.

For all $k\geq 1$, then,
$\Lk(K_{\gamma_{n_{k-1}}^{k-1}},K_\rho)=
p\,(A(k)q+a(k)s)$, and
Proposition~\ref{prop:first_general_slope_calculation}
gives the slope of $\gamma_{n_k}^k$ to
be $2p(A(k)q+a(k)s)+1/n_k$.

We now consider the second case $\tau \serho
\gamma_{n_0}^0 \netau \gamma_{n_1}^1 \netau \gamma_{n_2}^2 \netau \cdots$
of the drop-$\rho$ iteration. For the iterative step, when computing the
slope of $\gamma_{n_{k+1}}$,
the setup picture Figure~\ref{fig:first_general}(a) has $K_U=K_\tau$
and $K_L=K_{\gamma_{n_k}^k}$, the latter oriented so that its top
portion is $K_\tau$ oriented left-to-right, and bottom portion, originally
$K_{\gamma_{n_{k-1}}^{k-1}}$, has top portion (from $K_\tau$)
oriented left-to-right or right-to-left
according as $n_k$ is odd or even.
For $k=0$, $K_{\gamma_{n_0}}$ with this orientation represents
\[ (p+r,q+s)+\epsilon(0)(p,q) = (A(1)-a(1))(p+r,q+s) + a(1)(p,q)\ .\]
Inductively, assume that
$K_{\gamma_{n_{k-1}}^{k-1}}$ represents $(A(k)-a(k))(p+r,q+s)+a(k)(p,q)$.
Then, with the needed orientation for the setup picture,
$K_{\gamma_{n_k}^k}$ represents
\begin{gather*}
(p+r,q+s)+\epsilon(k)((A(k)-a(k))(p+r,q+s)+a(k)(p,q))\\
=(A(k+1)-a(k+1))(p+r,q+s)+a(k+1)(p,q)
\end{gather*}
The slope calculation of $\gamma_{n_k}^k$ is then
\begin{gather*} 2\Lk(K_\tau,K_{\gamma_{n_{k-1}}^k})+1/n_k
=2(q+s)((A(k)-a(k))(p+r)+a(k)p) + 1/n_k\\
=2p(q+s)A(k)+2r(q+s)(A(k)-a(k)) + 1/n_k\ .
\end{gather*}

These calculations have established the first two cases of the following
result. Each of the remaining six cases is very similar to one of the first
two. Summarizing, we have
\begin{theorem}\label{thm:iterated_slopes} The slopes
of the tunnels in the iterated splitting sequences for the torus knot
$T_{p+r,q+s}$ are as follows.
\vspace*{2 ex}
\begin{small}
\begin{center}
\renewcommand{\arraystretch}{1.3}
\setlength{\fboxsep}{0pt}
\setlength{\tabcolsep}{8pt}
\fbox{%
\begin{tabular}{c|c}
$\mathrm{sequence}$&$\mathrm{slope\ of}$ $\gamma_{n_k}^k$\\
\hline\hline
$\tau \serho \gamma_{n_0}^0 \serho \gamma_{n_1}^1 \serho \cdots$&%
$2p\,(\,A(k)q+a(k)s\,)+1/n_k$\\
$\tau \serho \gamma_{n_0}^0 \netau \gamma_{n_1}^1 \netau \cdots$&%
$2(q+s)\,(\,A(k)p+(A(k)-a(k))r\,) + 1/n_k$\\
$\tau \selambda \gamma_{n_0}^0 \selambda \gamma_{n_1}^1 \selambda \cdots$&%
$2r\,(\,A(k)s+a(k)q\,)+1/n_k$\\
$\tau \selambda \gamma_{n_0}^0 \netau \gamma_{n_1}^1 \netau \cdots$&%
$2(q+s)\,(\,A(k)r+(A(k)-a(k))p\,)+1/n_k$\\
$\tau \nerho \gamma_{n_0}^0 \nerho \gamma_{n_1}^1 \nerho \cdots$&%
$2q\,(\,A(k)p+a(k)r\,)+1/n_k$\\
$\tau \nerho \gamma_{n_0}^0 \setau \gamma_{n_1}^1 \setau \cdots$&%
$2(p+r)\,(\,A(k)q+(A(k)-a(k))s\,) + 1/n_k$\\
$\tau \nelambda \gamma_{n_0}^0 \nelambda \gamma_{n_1}^1 \nelambda \cdots$&%
$2s\,(\,A(k)r+a(k)p\,)+1/n_k$\\
$\tau \nelambda \gamma_{n_0}^0 \setau \gamma_{n_1}^1 \setau \cdots$&%
$2(p+r)\,(\,A(k)s+(A(k)-a(k))q\,)+1/n_k$\\
\end{tabular}}
\end{center}
\end{small}
\end{theorem}

The binary invariants produced by splitting and iterated splitting are
easily determined. When $\rho$ is one of the disks of the principal pair of
a tunnel (that is, one of the two disks in the principal vertex other than
the tunnel disk itself), a drop-$\rho$ or lift-$\rho$ splitting or
iterative step retains $\rho$ and replaces the other disk of the principal
pair. Thus, for example, in the all drop-$\rho$ iteration, every binary
invariant is $0$ except possible that of the splitting, which depends on
the cabling construction that preceded it (that is, the invariant is $0$ if
$\rho$ was in the principal pair of the tunnel for the cabling construction
that preceded the splitting, and $1$ if $\rho$ was the previous tunnel). In
a sequence such as $\tau \serho \gamma_{n_0}^0 \netau \gamma_{n_1}^1 \netau
\cdots$, the second binary invariant, associated to the first lift-$\tau$
step of the iteration, has binary invariant $1$, and all others
except possibly the initial splitting have binary invariant~$0$.

Since a splitting-and-iteration sequence can never have more than two
binary invariants equal to $1$, with the two $1$'s contiguous in that
case, the sequence can never increase the depth by more than $1$ from that
of the starting torus tunnel (see for example the last paragraph of
Section~3 of~\cite{CMbridge}). 

\section{Two-bridge knots}
\label{sec:2bridge_iteration}

A good example of the iterated splitting construction is furnished by
$2$-bridge knots. Indeed, in some sense the iterated splitting construction
is a far-reaching generalization of $2$-bridge knots. In this section, we
will see that any drop-$\rho$ iteration of the first kind examined in
Sections~\ref{sec:iterated} and~\ref{sec:iterated_slopes} and starting with
the trivial knot positioned as $T_{1,1}$ produces a $2$-bridge knot in the
$(1,1)$-position whose upper tunnel is the upper semisimple tunnel of the
knot, and moreover that every semisimple tunnel of every $2$-bridge
knot can be obtained in this way.

We wil use the notation and the description of the classification of
$2$-bridge knots presented in \cite[Section 10]{CMsemisimple}. We first
recall the calculation of the slope invariants of the upper semisimple
tunnel of a $2$-bridge knot given in~\cite[Proposition~10.4]{CMsemisimple}:
\begin{proposition}\label{prop:2bridge_slopes}
Let $K$ be a $2$-bridge knot in the $2$-bridge position corresponding to
the continued fraction $[2a_d,2b_d,\ldots,2a_0,2b_0]$, with $b_0\neq 0$ and
each $a_i=\pm 1$. Then the slope invariants of the upper semisimple tunnel
of $K$ are as follows:
\begin{enumerate}
\item[(i)] $m_0=\left[\displaystyle\frac{2b_0}{4b_0+1}\right]$ or
$m_0=\left[\displaystyle\frac{2b_0-1}{4b_0-1}\right]$
according as $a_0$ is $1$ or $-1$.
\item[(ii)] For $1\leq i\leq d$, $m_i=-2a_{i-1}+1/k_i$, where
\begin{enumerate}
\item[(a)] $k_i=2b_i+1$ if $a_i=a_{i-1}=1$,
\item[(b)] $k_i=2b_i$ if $a_i$ and $a_{i-1}$ have opposite signs, and
\item[(c)] $k_i=2b_i-1$ if $a_i=a_{i-1}=-1$.
\end{enumerate}
\end{enumerate}
\label{prop:semisimple_slopes}
\end{proposition}

Fix $K$ as in Proposition~\ref{prop:2bridge_slopes}. Denote the slope
invariants of its upper semisimple tunnel as given in
Proposition~\ref{prop:2bridge_slopes} by $m_0,\ldots\,$, $m_d$.

Starting with the trivial knot $T_{1,1}$, we will carry out a drop-$\rho$
splitting and iteration, that is, the first type detailed in each of
Sections~\ref{sec:iterated} and~\ref{sec:iterated_slopes}.  We have
\[M_{1,1}=\begin{pmatrix}1&0\\0&1\end{pmatrix} = I\ ,\]
thus $(p,q)=(1,0)$, $(r,s)=(0,1)$, and $K_\rho=T_{1,0}$.

Perform the initial drop-$\rho$ splitting with $n_0$ equal to $2b_0$ if
$a_0=1$ and to $2b_0-1$ if $a_0=-1$. Note that every nonzero choice of
$n_0$ occurs for some $m_0$. By Proposition~\ref{prop:sigma_slopes}(c) (or
Theorem~\ref{thm:iterated_slopes} with $k=0$), the slope of
$\gamma_{n_0}^0$ is $2+1/n_0$, so its simple slope is $[n_0/(2n_0+1)]$. By
Proposition~\ref{prop:2bridge_slopes}(i), this is $m_0$.

Now we carry out the first $d$ steps of the iteration, using $n_r=k_r$ at
each step. Again, every possible nonzero value of $n_r$ occurs for some
choice of $K$. We have $m_1=-2a_0+1/k_1$. If $a_0=1$, then $n_0$ was even
and (using the notation of Section~\ref{sec:iterated_slopes})
$a(1)=(-1)^{1+n_0}=-1$. If $a_0=-1$, then $n_0$ was odd and $a(1)=1$. In
either case, $a(1)=-a_0$. Theorem~\ref{thm:iterated_slopes} gives the slope
of $\gamma_{n_1}^1$ to be $2a(1)+1/n_1=-2a_0+1/k_1=m_1$.

For $r\geq 2$, assume inductively that $a(r)=-a_{r-1}$. If $n_r=k_r$ is even,
then we are in Case~(ii)(b) of
Proposition~\ref{prop:2bridge_slopes}, and $a_{r-1}=-a_r$. We find that
$a(r+1)=(-1)^{1+n_r}a(r)=-a(r)=a_{r-1}=-a_r$.
If $n_r$ is odd, then we are in Case~(ii)(a) or~(ii)(c) of
Proposition~\ref{prop:2bridge_slopes}, and $a_{r-1}=a_r$. We find that
$a(r+1)=(-1)^{1+n_r}a(r)=a(r)=-a_{r-1}=-a_r$, completing the induction.

Theorem~\ref{thm:iterated_slopes} now gives the slope of $\gamma_{n_r}^r$
to be \[ 2a(r)+1/n_r = -2a_{r-1}+1/k_r=m_r\ ,\] completing the induction.

\section{Classification of tunnels}
\label{sec:TC}

At this point in history one may begin to contemplate a classification of
tunnels of tunnel number 1 knots based on the examples that have been found
during the past several decades. In this section we will list the cases
that occur or appear to occur. It is plausible that this list may be
complete or nearly so, but we are unaware of any evidence supporting this
other than the absence of other examples found and a sense that there ought
to be a fairly strict limitation on the complexity of tunnel behavior for a
given knot.

In our list it is to be understood that in some cases, tunnels are
equivalent due to symmetries or degeneracies. For example, the upper and
lower tunnels of a $2$-bridge knot may be equivalent under an involution of
$S^3$ preserving the knot, and the middle tunnel of a torus knot is known
to be isotopic to the upper or lower tunnel for certain cases (and hence is
a $(1,1)$-tunnel rather than a regular tunnel).

We list the cases, then comment on them below.
\begin{KTP} These are the known possibilities for the set of tunnels of
a tunnel number $1$ knot $K$, allowing some of the tunnels to be
equivalent due to symmetries or degeneracies:
\begin{enumerate}
\item[I.] $K$ has a unique regular tunnel.
\item[II.] $K$ has one $(1,1)$-position and two $(1,1)$-tunnels.
\item[III.] $K$ has two $(1,1)$-positions and four $(1,1)$-tunnels.
\item[IV.] $K$ has one $(1,1)$-position and two $(1,1)$-tunnels, plus one
regular tunnel.
\item[V.] $K$ has two $(1,1)$-positions and four $(1,1)$-tunnels, plus one
regular tunnel.
\item[VI.] $K$ has one $(1,1)$-position and two $(1,1)$-tunnels, plus two
regular tunnels.
\item[VII.] $K$ has no $(1,1)$-position, but has two regular tunnels.
\end{enumerate}
\end{KTP}

We now comment on the individual cases.
\smallskip

\noindent \textsl{Case I} 

As explained in \cite[Section~3]{CMbridge}, results of M. Scharlemann and
M. Tomova~\cite{Scharlemann-Tomova} and J. Johnson~\cite{Johnson} combine
to show that whenever $K$ has a tunnel of Hempel distance at least $6$
(that is, the Hempel distance of the associated genus-$2$ Heegaard
splitting of the exterior of $K$), it is the unique tunnel of $K$. Thus
Case~I holds for all high-distance tunnels.
\smallskip

\newpage
\noindent \textsl{Case II} 

This seems likely to be the generic case when $K$ has a $(1,1)$-tunnel,
although we are not aware of any examples for which it has been proven that
a specific knot admits exactly two $(1,1)$-tunnels, other than symmetric or
degenerate cases such as torus knots for which the middle tunnel is
equivalent to the upper or lower $(1,1)$-tunnel.
\smallskip

\noindent \textsl{Case III} 

Tunnels of $2$-bridge knots are fully classified due to work of several
authors, and they satisfy Case~III. D. Heath and H. Song~\cite{Heath-Song}
proved that the $(-2,3,7)$-pretzel knot satisfies Case~III, and there are
expected to be other examples.
\smallskip

\noindent \textsl{Case IV} 

Torus knots and their middle tunnels are the long-known examples of
Case~IV. Assuming that at least some of them have no other unknown tunnels,
the examples generated in \cite{CMsplitting} and this paper provide more
such knots. See also the comments on the remaining three cases.\par
\smallskip

\noindent \textsl{Cases V, VI, and VII} 

These remaining cases describe examples recently found and kindly provided
to us by John Berge~\cite{Berge}. They were obtained using his software
\textit{Heegaard,} which works with two-generator one-relator presentations
of $\pi_1(S^3-K)$ whose generators are free generators of the fundamental
group of the exterior handlebody $H'=\overline{S^3-H}$, and whose relator
is represented by the boundary $C$ of a tunnel disk $D$ in $H$. The knot
$K$ is the usual knot associated to $D$, that is, a core circle of the
solid torus $\overline{H- N(D)}$, where $N(D)$ is a regular neighborhood of
$D$ in $H$. \textit{Heegaard} is able to distinguish equivalence classes of
such $C$ under diffeomorphism of $H'$, showing that the tunnel disks they
bound cannot be equivalent. Regularity of the tunnels can be tested by
using a procedure (also used by K. Ishihara~\cite{Ishihara}) that finds the
principal meridian pair for $K$ associated to a tunnel, and then checking
whether either of the disks is primitive; primitivity of a disk $E\subset
H$ in our sense (that is, $\partial E$ crosses the boundary of some disk
$E'\subset H'$ exactly once) is equivalent to primitivity of $\partial E$
as an element of $\pi_1(H')$, and can be checked algebraically.

Once a tunnel has been found, the software searches for more tunnels for
the knot by a method that generates a large number of additional such
two-generator one-relator presentations for $\pi_1(S^3-K)$ and tests them
for isomorphism with those already found. Although there is no known means
to ensure that this method finds all of the tunnels for these examples, it
seems likely that it does. For example, for the $(-2,3,7)$-pretzel knot,
all four tunnels are found among the first few of the large number of
presentations that the software examines.

Berge examined the hyperbolic double-primitive knots $K$ having Dehn
surgeries that produce lens spaces $L(p,q)$ with $p<100$, and the
``sporadic'' double-primitive knots of Types 9, 10, 11, and 12 (detailed in
J. Berge~\cite{doubleprimitive}) having Dehn surgeries that produce lens
spaces $L(p,q)$ with $p<500$, as well as some non-double-primitive
knots. Assuming that the software did find all tunnels of those knots, the
possibilities listed in Cases~V, VI, VII were obtained, as well as quite a
few instances of the other cases including Case~IV. Some of the examples of
Case~VII occurred for knots that are not double-primitive. We do not know
whether the regular tunnels in his examples of Cases~IV, V, and~VI arise
from $(1,1)$-positions by the construction we have examined in this paper.

\bibliographystyle{amsplain}

\begin{thebibliography}{10}

\bibitem{Berge} J. Berge, personal communication. He may be contacted at
jberge@charter.net for additional information about \textit{Heegaard.}

\bibitem{doubleprimitive} J. Berge, Some knots with surgeries yielding lens
spaces, preprint.

\bibitem{B-R-Z} M. Boileau, M. Rost, and H. Zieschang, On Heegaard
decompositions of torus knot exteriors and related Seifert fibre spaces,
\textit{Math. Ann.} 279 (1988), 553--581.

\bibitem{CMtree} S. Cho and D. McCullough, The tree of knot tunnels,
\textit{Geom. Topol.} 13 (2009) 769-815.

\bibitem{CMtorus} S. Cho and D. McCullough, Cabling sequences of
tunnels of torus knots, \textit{Algebr. Geom. Topol.} 9 (2009) 1--20.

\bibitem{CMgiant_steps} S. Cho and D. McCullough, Constructing knot tunnels
using giant steps, \textit{Proc. Amer. Math. Soc.} 138 (2010), 375-384.

\bibitem{CMbridge} S. Cho and D. McCullough, Tunnel leveling, depth, and
bridge numbers, \textit{Trans. Amer. Math. Soc.} 353 (2011), 259--280.

\bibitem{CMsemisimple} S. Cho and D. McCullough, Semisimple tunnels,
arXiv:1006.5232.

\bibitem{CMsplitting} S. Cho and D. McCullough, Middle tunnels by
splitting, arXiv:1108.3425

\bibitem{CMsoftware} S. Cho and D. McCullough, software available at
\texttt{math.ou.edu/$_{\widetilde{\phantom{i}}}\,$dmccullough}~.

\bibitem{Goda-Hayashi} H. Goda and C. Hayashi, Genus two Heegaard
splittings of exteriors of 1-genus 1-bridge knots, to appear in
\textit{Kobe J. Math.}

\bibitem{Heath-Song} D. Heath and H.-J.\ Song, Unknotting tunnels for $P(-2,3,7)$,
\textit{J. Knot Theory Ramifications} 14 (2005), 1077--1085.

\bibitem{Ishihara} K. Ishihara, An algorithm for finding parameters of
tunnels, \textit{Alg.\ Geom.\ Topology.} 11 (2011), 2167--2190.

\bibitem{Johnson} J. Johnson, Bridge number and the curve
complex, arXiv math.GT/\allowbreak0603102.

\bibitem{Moriah} Y. Moriah, Heegaard splittings of Seifert fibered spaces,
\textit{Invent. Math.} 91 (1988), 465--481.

\bibitem{Morimoto-Sakuma} K. Morimoto and M. Sakuma, On unknotting
tunnels for knots, \textit{Math. Ann.} 289 (1991), 143--167.

\bibitem{Morimoto-Sakuma-Yokota} K. Morimoto, M. Sakuma, and
Y.Yokota, Examples of tunnel number $1$ knots which have the ``$1+1=3$''
property, \textit{Math. Proc. Camb. Phil. Soc.} 119 (1996), 113--118.

\bibitem{Scharlemann-Tomova} M. Scharlemann and M. Tomova, Alternate
Heegaard genus bounds distance, \textit{Geom. Topol.} 10 (2006), 593--617.

\end{thebibliography}

\end{document}